\documentclass{cmbe17-bis}

\usepackage{amsmath,amsfonts,amssymb}
\usepackage{graphicx,wrapfig}

\usepackage[hyphens]{url}
\usepackage[colorlinks]{hyperref}
\usepackage{subcaption}
\usepackage{pgfplotstable}
\usepackage{xspace}
\usepackage{bm}
\usepackage{todonotes}
\usepackage{pgfplots}
\pgfqkeys{/pgfplots}{
  legend style = { at={(0.95,0.3)} },
  grid=major
}
\usepackage{pgf}
 \pgfplotstableset{
 columns/{tools}/.style={string type,column name=\textsc{Tools},column type=l},
  columns/{Li}/.style={string type,column name=\textsc{$L_i$},column type=l},
  columns/{Lo}/.style={string type,column name=\textsc{$L_o$},column type=l},
    columns/{Time}/.style={string type,column name=\textsc{Time},column type=l},
    columns/{Tf}/.style={string type,column name=\textsc{$T_f$},column type=l},
    columns/{Machine}/.style={string type,column name=\textsc{Machine},column type=c},
    columns/{Reynolds}/.style={string type,column name=\textsc{$\re$},column type=l},
 columns/{geo}/.style={string type,column name=\textsc{Geometry},column type=l},
columns/{tetras}/.style={column name=\textsc{Tetrahedrons},column type=c,int detect},
columns/{NDOF}/.style={column name=\textsc{DOF},column type=c,int detect},
columns/{NDOFL}/.style={column name=\textsc{$DOF_{local}$},column type=c,int detect},
columns/{NNZ}/.style={column name=\textsc{NNZ},column type=c,int detect},
columns/{Nproc}/.style={column name=\textsc{$N^{core}$},column type=c,int detect},
columns/{CPU}/.style={column name=\textsc{CPU},column type=c,int detect},
columns/{nodes}/.style={column name=\textsc{Nodes},column type=c,int detect},
columns/{triangles}/.style={column name=\textsc{Triangles},column type=c,int detect},
columns/{hsize}/.style={column name=h, string type, column type={l}},
columns/{hmin}/.style={column name=\textsc{$h_{min}$},column type=c,int detect},
columns/{hmax}/.style={column name=\textsc{$h_{max}$},column type=c,int detect},
columns/{havrg}/.style={column name=\textsc{$h_{avrg}$},column type=c,int detect},
columns/{ReThroat}/.style={column name=\textsc{$\re_{t}$},column type=l,int detect},
columns/{ReInlet}/.style={column name=\textsc{$\re_{i}$},column type=l,int detect},
columns/{flowrate}/.style={column name=\textsc{Flow rate},column type=c,int detect},
columns/{TimeStep}/.style={column name=$\Delta t$},
columns/{errorUL2}/.style={column name=$\|\bu-\bu_\delta\|_{0,\Omega_\delta}$,sci,sci zerofill,precision=1},
columns/{errorUH1}/.style={column name=$\|\bu-\bu_\delta\|_{1,\Omega_\delta}$,sci,sci zerofill,precision=1},
columns/{errorP}/.style={column name=$\|p-p_\delta\|_{0,\Omega_\delta}$,sci,sci zerofill,precision=1},
columns/{F1In}/.style={column name={$||{\bm F}_{ex}-{\bm F}_\delta||$},sci,sci zerofill,precision=1},
create on use/Fin/.style={create col/expr={sqrt{\thisrow{F1In}^2}}},
empty cells with={--},
string replace={nan}{}, 
string replace={-}{-}, 
every head row/.style={before row=\toprule,after row=\midrule},
every last row/.style={after row=\bottomrule}}
\usepackage{booktabs}

\newcommand{\bu}{\mathbf{u}}

\newcommand{\kgeo}{\ensuremath{k_\mathrm{geo}}\xspace}

\newcommand{\re}{\mathcal{R}e}

\newcommand{\ppg}[3]{\ensuremath{P_{#1}P_{#2}G_{#3}}\xspace}

\title{High order finite element simulations for fluid dynamics validated by experimental data from the FDA benchmark nozzle model}

\author{V. Chabannes}
\author{C. Prud'homme}
\author{M. Szopos}
\author{R. Tarabay}
\affil{Universit\'e de Strasbourg, CNRS, IRMA UMR 7501, F-67000 Strasbourg, France, \texttt{\{chabannes,prudhomme,szopos,tarabay\}@math.unistra.fr}}

\summary{The objective of the present work is to construct a sound mathematical, numerical and computational framework relevant to blood flow simulations and
to assess it through a careful validation against experimental data. We perform simulations of a benchmark proposed by the FDA
for fluid flow in an idealized medical device, under different flow regimes. The results are evaluated using metrics proposed in the
literature and the findings are in very good agreement with the validation experiment.}


\begin{document}

\section{Introduction}


\begin{wrapfigure}{r}{.4\linewidth}
  \begin{center}
    \begin{subfigure}[b]{.97\linewidth}
      \includegraphics[width=0.8\textwidth]{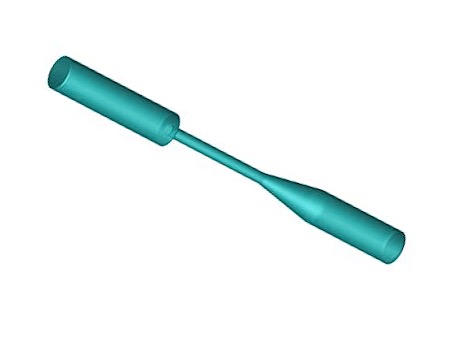}
    \end{subfigure}\\
    \begin{subfigure}[b]{.97\linewidth}
      \includegraphics[width=0.8\textwidth]{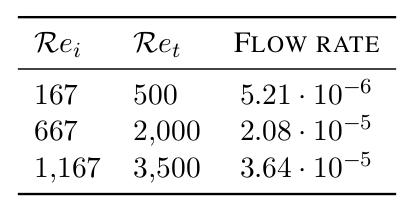}
   \end{subfigure}
   \caption{Computational domain (top) and flow regime specifications (bottom); $\re_i$ and $\re_t$: Reynolds number in the inlet section and throat section, respectively.}
   \label{fig:comp_domain_Re}
  \end{center}
\end{wrapfigure}

\noindent A challenging benchmark was proposed by the US Food and Drug
Administration (FDA) in \cite{hariharan2011} in order to assess the
stability, accuracy and robustness of computational methods in
different physiological regimes. The findings of 28 blinded
investigations were reported in \cite{stewart2012} and as critically
analyzed in \cite{sotiropoulos2012}, practically all CFD solvers
failed to predict results that agreed in a satisfactory manner with
the experimental data. Several subsequent papers tackled this
question, by employing different numerical approaches: for instance a
finite-element based direct numerical simulation method in
\cite{passerini2013} or a large-eddy simulation method in
\cite{zmijanovic2016}.

\noindent We aim at contributing to the effort of improving the
reliability and reproducibility of computational studies by performing
a thorough validation of the fluid solver developed in the open source
finite element library Feel++~\cite{feelpp}. In the current
investigation, we present results corresponding to three Reynolds
numbers $500,2000$ and $3500$ obtained by using a direct numerical
simulation method of the Navier-Stokes equations. In particular we implement and compare low order
as well as high order approximations including for the geometry and we
discuss some issues not previously reported in the literature.


%
%

\section{Methodology}

\noindent{\bf Benchmark description.} The FDA benchmark nozzle model
provides a comprehensive dataset of experimental measures using a
well-defined geometry corresponding to an idealized medical device
(see Figure~\ref{fig:comp_domain_Re} for a schematic sketch of the
domain and \cite[Sec. 2.1]{hariharan2011} for the precise dimensions
of each part). Five sets of data spanning laminar, transitional and
turbulent regimes are made available; we focus in the current work on
the flow regime specifications described in
Figure~\ref{fig:comp_domain_Re}.


\noindent The comparison with experimental data is made in terms of (i) wall pressure difference (normalized to the mean throat velocity) versus axial distance;
 and (ii) axial component of the velocity (normalized to the mean inlet velocity) along the centerline:\\
\begin{equation}
 \displaystyle \Delta p^{norm} = \frac{p_z - p_{z=0}}{\frac{1}{2} \rho_f \overline{u}_t^2} \text{ and } \quad
 u_z^{norm} = \frac{u_z}{\overline{u}_i},\text{ where }\quad \displaystyle \overline{u}_t^2 = \frac{4Q}{\pi D_t^2},\ \displaystyle \overline{u}_i = \frac{4Q}{\pi D_i^2},
 \label{eq:press_vel_norm}
\end{equation}
and $Q$ is the volumetric flow rate retrieved from $\re_t$ (see
Figure~\ref{fig:comp_domain_Re}, right panel). Furthermore, we present
results on two validation metrics reported in \cite{stewart2012}, also
assessed in \cite{passerini2013}: a conservation of mass error metric
$E_Q$ (on a percentage basis) and a general validation metric $E_z$
comparing average experimental velocity data with computed axial
velocities.
%
%
%

{\bf Fluid equations and numerical approach.}
We now turn to the mathematical and the numerical setting. We consider
the homogeneous, incompressible, unsteady Navier-Stokes equations,
which read in conservative form: find $(\mathbf{u},p)$ such that
$\rho \left( \frac{\partial \mathbf{u}}{\partial t} + (\mathbf{u} \cdot
  \nabla) \mathbf{u} \right) - \mu \Delta \mathbf{u} + \nabla p     =  \mathbf{0}, \quad \text{div} (\mathbf{u}) =  0,$
in $\Omega \times I$. The set $\Omega \subset \mathbb{R}^3$ represents
the spatial domain described in Figure~\ref{fig:comp_domain_Re}.
{$I=\left(0,T\right)$} is the time interval, $\mathbf{u}$ and $p$ are
the velocity and pressure of the fluid
and $\rho$ and $\mu$ are the
density and the dynamic viscosity of the fluid, respectively.  We
supplement the equations with initial and boundary conditions. At
$t=0s$, the fluid is considered to be at rest,
$\mathbf{u}(\mathbf{x},t)=\mathbf{0}$. A Poiseuille velocity profile is imposed on
$\Gamma_{\mathrm{inlet}}$, homogeneous Dirichlet condition on $\Gamma_{\mathrm{wall}}$ and a free outflow on
$\Gamma_{\mathrm{outlet}}$.

We refer to \cite[Sec. 2]{vivabrain-cemracs13} regarding the variational
formulation, the finite element discretization including low to high
order geometry as well as the time discretization. We choose the generalized Taylor-Hood finite
element for the velocity-pressure discretization; the notation \ppg{N+1}{N}{\kgeo} is used to specify
exactly the discretization spaces for the velocity, pressure and geometry, respectively.

The benchmark hereafter is developed in the framework of the Finite
Element Embedded Library in C++, Feel++\cite{feelpp}, that allows to
use a very wide range of Galerkin methods and advanced numerical
techniques such as domain decomposition. The ingredients include a
very expressive embedded language, seamless interpolation, mesh
adaption and seamless parallelization. Regarding the computational
domain, we used \textsc{Gmsh}. The construction used the following
steps: \textit{(i)} start with a 2D geometry embedding the benchmark
metric locations and customizing characteristic mesh size depending on
the region and \textit{(ii)} extrude by rotation to obtain the device
geometry.  Finally we use the \textsc{PETSc} interface developed in
Feel++ and in particular the FieldSplit preconditioning framework to
implement block preconditioning strategies such as
PCD~\cite{elman2014}. Note that PCD requires specific tuning with
respect to boundary conditions.

\begin{wraptable}{r}{.6\linewidth}
  \centering
 \begin{tabular}{crrrrr}
  \toprule
    & $h_{min}$ & $h_{max}$ & $h_{average}$ & $N_{\mathrm{elt}}$ \\
  \midrule
  \textsc{m0} & $1.9\cdot 10^{-4}$ & $2.9\cdot 10^{-3}$ & $1.3\cdot 10^{-3}$ & $412\,575$   \\
  \textsc{m1} & $1.6\cdot 10^{-4}$ & $1.8\cdot 10^{-3}$ & $7.6\cdot 10^{-4}$ & $830\,000$   \\
  \textsc{m2} & $1.4\cdot 10^{-4}$ & $1.96\cdot 10^{-3}$ & $6\cdot 10^{-4}$ & $3\,400\,000$   \\
  \textsc{m3} & $8.5\cdot 10^{-5}$ & $1.7\cdot 10^{-3}$ & $3.5\cdot 10^{-4}$ & $7\,000\,000$   \\
  \textsc{m4} & $6.3\cdot 10^{-5}$ & $2.0\cdot 10^{-3}$ & $5.8\cdot 10^{-4}$ & $2\,879\,365$   \\
  \textsc{m5} & $1.4\cdot 10^{-4}$ & $2.6\cdot 10^{-3}$ & $4.1\cdot 10^{-4}$ & $3\,200\,000$   \\
\bottomrule
 \end{tabular}
\caption{Characteristic lengths of the different meshes: $h_{min}, h_{max},
  h_{average}$ are respectively the minimum, maximum and average edge
  length in the meshes and $N_{\mathrm{elt}}$ is the number
  of tetrahedra.}
\label{tab:results:mesh:characteristic}
\end{wraptable}

\section{Results and conclusions}
\label{sec:results-conclusions}



We perform simulations for three Reynolds numbers evaluated in the throat $\re_t=500,2000,3500$, with several mesh refinements
and polynomial order approximations. The fluid's prescribed density is $\rho = 1056 \frac{kg}{m^3}$ and viscosity $\mu = 0.0035 Pa.s$.
The mesh characteristics are described in Table~\ref{tab:results:mesh:characteristic}.
At $\re_t=500$, the simulation is carried out until $t=3 s$, time reasonably close to the steady state,
and we choose the time step equal to $\Delta t=10^{-3}$.
At $\re_t=2000$ (resp $\re_t=3500$), the numerical experiments were carried out until $t=0.45 s$ (resp $t=0.4 s$),
time when the turbulent regime was fully developed and we set $\Delta t=10^{-4}$.

Figure~\ref{fig:results:velocity3500} shows the results in the three flow
regimes for the normalized axial velocity and the normalized pressure difference along the $z$ axis,
respectively. In each case, we can see satisfactory agreement with the experimental data. However,
for $\re_t=2000$, we observe that the numerical jet breakdown point is
captured further downstream than the experimentally observed
breakdown point. As recently highlighted in \cite{zmijanovic2016}, the prediction of the axial location of
the jet breakdown is extremely sensitive to numerical parameters, therefore a possible explanation of this
mismatch may be the accuracy of the numerical integration. Finally, we illustrate in Figure~\ref{fig::results:metrics500} the
computation of metrics $E_z$ and $E_Q$ for several mesh refinements at
$\re_t=500$.  The metric $E_z$ takes small values in each numerical experiment, identifying a
good agreement betweeen computed and experimental data, and
displays only small variations with respect to mesh refinement. On the other hand, the metric $E_Q$
is more sensitive to this factor: error doesn't exceed the $\sim 2\%$ except for the coarse mesh $M0$
where, in two locations, the error increases up to $\sim10\%$.  Furthermore,
we note that the \ppg{3}{2}{1} approximation doesn't improve the results for the coarse mesh, but that a satisfactory
error below $ 2\%$ is retrieved when using a \ppg{2}{1}{2} approximation. Additional tests to complement the
study of the impact of high order approximation are ongoing.

{\bf Conclusions and perspectives} We validated our
computation fluid dynamic framework against this FDA benchmark for
three different regimes and different discretization and solution
strategies. Perspectives include a full report on our findings
including in terms of iteration and timing performances as well
extending our results to the turbulent range.

{\bf Acknowledgments} The authors wish to thank Mourad Ismail for the
fruitful discussions we had. Moreover we would like to acknowledge the support of
\textit{(i)} Center of Modeling and Simulation of Strasbourg
(Cemosis), \textit{(ii)} the ANR MONU-Vivabrain \textit{(iii)} the
LabEx IRMIA and  \textit{(iv)} PRACE for awarding us access to resource
Curie based in France at CCRT as well as GENCI for awarding us access
to resource Occigen based in France at Cines

\begin{figure}[htbp]
    \begin{subfigure}[b]{.47\linewidth}
      \begin{tikzpicture}[scale=0.85]
        \begin{axis}[
            xmin=-0.1,xmax=0.1,
            xlabel=arc length,
            ylabel=normalized $u_z$,
            legend style = {legend pos = south east},
            legend columns=2
          ]
          \addplot table[x=arcLength,y expr=\thisrow{normalizedUz}/0.04603291471] {Figures/dats/tarabay/FDA/ExpData/Re500/Sudden_Expansion_500_243_velocity.dat};
          \addplot table[x=arcLength,y expr=\thisrow{normalizedUz}/0.04603291471]{Figures/dats/tarabay/FDA/ExpData/Re500/Sudden_Expansion_500_297_velocity.dat};
          \addplot table[x=arcLength,y expr=\thisrow{normalizedUz}/0.04603291471]{Figures/dats/tarabay/FDA/ExpData/Re500/Sudden_Expansion_500_468_velocity.dat};
          \addplot table[x=arcLength,y expr=\thisrow{normalizedUz}/0.04603291471]{Figures/dats/tarabay/FDA/ExpData/Re500/Sudden_Expansion_500_763_velocity.dat};
          \addplot table[x=arcLength,y expr=\thisrow{normalizedUz}/0.04603291471]{Figures/dats/tarabay/FDA/ExpData/Re500/Sudden_Expansion_500_999_velocity.dat};
          \addplot[color=orange,dashed,line width=2pt] table[x=Points2,y expr=\thisrow{u2}/0.04603291471, col sep=comma]{Figures/dats/tarabay/FDA/Re500_Coarser_P2P1G2/P2P1G1_L_extrude_it3000_LU_pressure_z.csv};
          \addplot[color=green, dashed,line width=2pt] table[x=Points2,y expr=\thisrow{u2}/0.04603291471]{Figures/dats/tarabay/FDA/Re500_Fine_P2P1G1/P2P1G1_L_extrude_it3000_LU_pressure_z.csv};

          \legend{data1, data2, data3, data4,data5, \textsc{M0}-P2P1G2, \textsc{M3}-P2P1G1 }
        \end{axis}
      \end{tikzpicture}
    \end{subfigure}
    \begin{subfigure}[b]{.47\linewidth}
      \begin{tikzpicture}[scale=0.85]
        \begin{axis}[
            xmin=-0.1,xmax=0.1,
            ylabel near ticks,
            yticklabel pos=right,
            xlabel=arc length,
            ylabel=normalized pressure difference,
            legend style = {legend pos = north east},
            legend columns=2
          ]

          \addplot table[x=arcLength,y expr=(\thisrow{normalizedP}+130.879)/90.62664235] {Figures/dats/tarabay/FDA/ExpData/Re500/Sudden_Expansion_500_243_pressure.dat};
          \addplot table[x=arcLength,y expr=(\thisrow{normalizedP}-174.842)/90.62664235] {Figures/dats/tarabay/FDA/ExpData/Re500/Sudden_Expansion_500_297_pressure.dat};
          \addplot table[x=arcLength,y expr=(\thisrow{normalizedP}+144.963)/90.62664235] {Figures/dats/tarabay/FDA/ExpData/Re500/Sudden_Expansion_500_468_pressure.dat};
          \addplot table[x=arcLength,y expr=(\thisrow{normalizedP}+137.676)/90.62664235] {Figures/dats/tarabay/FDA/ExpData/Re500/Sudden_Expansion_500_763_pressure.dat};

          \addplot table[x=arcLength,y
          expr=(\thisrow{normalizedP}+7.26346)/90.62664235]
          {Figures/dats/tarabay/FDA/ExpData/Re500/Sudden_Expansion_500_999_pressure.dat};

          \addplot[color=orange,dashed,line width=2pt] table[x=Points2,y expr=(\thisrow{p}+44.277)/90.62664235, col sep=comma]{Figures/dats/tarabay/FDA/Re500_Coarser_P2P1G2/P2P1G1_L_extrude_it3000_LU_pressure_z.csv};
          \addplot[color=green, dashed,line width=2pt] table[x=Points2,y expr=(\thisrow{p}+28.406)/90.62664235]{Figures/dats/tarabay/FDA/Re500_Fine_P2P1G1/P2P1G1_L_extrude_it3000_LU_pressure_z.csv};
          \legend{data1, data2, data3, data4,data5,\textsc{M0}-P2P1G2,
            \textsc{M3}-P2P1G1}
        \end{axis}
      \end{tikzpicture}
    \end{subfigure}
    \begin{subfigure}[b]{.47\linewidth}
      \begin{tikzpicture}[scale=0.85]
        \begin{axis}[
            xmin=-0.1,xmax=0.1,
            xlabel=arc length,
            ylabel=normalized $u_z$,
            legend style={at={(0.7,0.)}, anchor=south east}
          ]
          \addplot table[x=arcLength,y expr=\thisrow{normalizedUz}/0.1842250631] {Figures/dats/tarabay/FDA/ExpData/Re2000/Sudden_Expansion_2000_243_velocity.dat};
          \addplot table[x=arcLength,y expr=\thisrow{normalizedUz}/0.1842250631]{Figures/dats/tarabay/FDA/ExpData/Re2000/Sudden_Expansion_2000_297_velocity.dat};
          \addplot table[x=arcLength,y expr=\thisrow{normalizedUz}/0.1842250631]{Figures/dats/tarabay/FDA/ExpData/Re2000/Sudden_Expansion_2000_468_velocity.dat};
          \addplot table[x=arcLength,y expr=\thisrow{normalizedUz}/0.1842250631]{Figures/dats/tarabay/FDA/ExpData/Re2000/Sudden_Expansion_2000_763_velocity.dat};
          \addplot table[x=arcLength,y expr=\thisrow{normalizedUz}/0.1842250631]{Figures/dats/tarabay/FDA/ExpData/Re2000/Sudden_Expansion_2000_999_velocity.dat};
          \addplot[green,dashed,line width=2pt] table[x=Points2,y expr=\thisrow{u2}/0.1842250631, col sep=comma]{Figures/dats/tarabay/FDA/latestFDA2000/P2P1G1_L_extrude_it3000_LU_pressure_z.csv};
          \legend{data1, data2, data3, data4,data5,\textsc{M4}-P2P1G1}
        \end{axis}
      \end{tikzpicture}
    \end{subfigure}
    \begin{subfigure}[b]{.47\linewidth}
      \begin{tikzpicture}[scale=0.85]
        \begin{axis}[
            xmin=-0.1,xmax=0.1,
            ylabel near ticks,
            yticklabel pos=right,
            xlabel=arc length,
            ylabel=normalized pressure difference,
            legend style = {legend pos = north east}
          ]

          \addplot table[x=arcLength,y expr=(\thisrow{normalizedP})/1451.497758] {Figures/dats/tarabay/FDA/ExpData/Re2000/Sudden_Expansion_2000_243_pressure.dat};
          \addplot table[x=arcLength,y expr=(\thisrow{normalizedP})/1451.497758] {Figures/dats/tarabay/FDA/ExpData/Re2000/Sudden_Expansion_2000_297_pressure.dat};
          \addplot table[x=arcLength,y expr=(\thisrow{normalizedP})/1451.497758] {Figures/dats/tarabay/FDA/ExpData/Re2000/Sudden_Expansion_2000_468_pressure.dat};
          \addplot table[x=arcLength,y expr=(\thisrow{normalizedP})/1451.497758] {Figures/dats/tarabay/FDA/ExpData/Re2000/Sudden_Expansion_2000_763_pressure.dat};
          \addplot[color=green,dashed,line width=2pt] table[x=Points2,y expr=(\thisrow{p}+337.11)/1451.497758, col sep=comma]{Figures/dats/tarabay/FDA/latestFDA2000/P2P1G1_L_extrude_it3000_LU_pressure_z.csv};
          \legend{data1, data2, data3, data4,\textsc{M4}-P2P1G1}
        \end{axis}
      \end{tikzpicture}
    \end{subfigure}\\
    \begin{subfigure}[b]{.47\linewidth}
      \begin{tikzpicture}[scale=0.85]
        \begin{axis}[
            xmin=-0.1,xmax=0.1,
            xlabel=arc length,
            ylabel=normalized $u_z$,
            legend style={at={(0.7,0.)}, anchor=south east}
          ]
          \addplot table[x=arcLength,y expr=\thisrow{normalizedUz}/0.322200348259] {Figures/dats/tarabay/FDA/ExpData/Re3500/Sudden_Expansion_500_243_velocity.dat};
          \addplot table[x=arcLength,y expr=\thisrow{normalizedUz}/0.322200348259]{Figures/dats/tarabay/FDA/ExpData/Re3500/Sudden_Expansion_500_297_velocity.dat};
          \addplot table[x=arcLength,y expr=\thisrow{normalizedUz}/0.322200348259]{Figures/dats/tarabay/FDA/ExpData/Re3500/Sudden_Expansion_500_468_velocity.dat};
          \addplot table[x=arcLength,y expr=\thisrow{normalizedUz}/0.322200348259]{Figures/dats/tarabay/FDA/ExpData/Re3500/Sudden_Expansion_500_763_velocity.dat};
          \addplot table[x=arcLength,y expr=\thisrow{normalizedUz}/0.322200348259]{Figures/dats/tarabay/FDA/ExpData/Re3500/Sudden_Expansion_500_999_velocity.dat};
          \addplot[color=green,dashed,line width=2pt] table[x=Points2,y expr=\thisrow{u2}/0.322200348259, col sep=comma]{Figures/dats/tarabay/FDA/Re3500_Large_P2P1G1/P2P1G1_L_extrude_it3000_LU_pressure_z.csv};
          \legend{data1, data2, data3, data4,data5,\textsc{M5}-P2P1G1}
        \end{axis}
      \end{tikzpicture}
    \end{subfigure}
    \begin{subfigure}[b]{.47\linewidth}
      \begin{tikzpicture}[scale=0.85]
        \begin{axis}[
            xmin=-0.1,xmax=0.1,
            ylabel near ticks,
            yticklabel pos=right,
            xlabel=arc length,
            ylabel=normalized pressure difference,
            legend style = {legend pos = north east}
          ]

          \addplot table[x=arcLength,y expr=(\thisrow{normalizedP})/4443.307946] {Figures/dats/tarabay/FDA/ExpData/Re3500/Sudden_Expansion_500_243_pressure.dat};
          \addplot table[x=arcLength,y expr=(\thisrow{normalizedP} -539.428)/4443.307946] {Figures/dats/tarabay/FDA/ExpData/Re3500/Sudden_Expansion_500_297_pressure.dat};
          \addplot table[x=arcLength,y expr=(\thisrow{normalizedP}-134.823)/4443.307946] {Figures/dats/tarabay/FDA/ExpData/Re3500/Sudden_Expansion_500_468_pressure.dat};
          \addplot table[x=arcLength,y expr=(\thisrow{normalizedP})/4443.307946] {Figures/dats/tarabay/FDA/ExpData/Re3500/Sudden_Expansion_500_763_pressure.dat};
          \addplot table[x=arcLength,y expr=(\thisrow{normalizedP}+9.05567)/4443.307946] {Figures/dats/tarabay/FDA/ExpData/Re3500/Sudden_Expansion_500_999_pressure.dat};
          \addplot[color=green,dashed,line width=2pt] table[x=Points2,y expr=(\thisrow{p}+1314.82)/4443.307946, col sep=comma]{Figures/dats/tarabay/FDA/Re3500_Large_P2P1G1/P2P1G1_L_extrude_it3000_LU_pressure_z.csv};
          \legend{data1, data2, data3, data4,data5,\textsc{M5}-P2P1G1}
        \end{axis}
      \end{tikzpicture}
    \end{subfigure}
    \caption{Comparison between experimental data and numerical
      results for the normalized axial velocity along z (left) and the
      normalized pressure difference along z (right), for $\re_t=500$
      (top) $\re_t=2000$ (middle) and $\re_t=3500$ (bottom).}
  \label{fig:results:velocity3500}

    \begin{subfigure}[b]{.47\linewidth}
      \begin{tikzpicture}[scale=0.85]

        \begin{semilogyaxis}[
            ymax=10,
            xlabel=z,
            ylabel=$E_z$,
            legend style = {legend pos = north west}
            ]
          \addplot[color=green,mark=square*] table[x=h,y=EzCP2P1]{Figures/dats/tarabay/FDA/Ez.dat};
          \addplot[color=orange,mark=square*] table[x=h,y=EzCP2P1G2] {Figures/dats/tarabay/FDA/Ez.dat};
          \addplot[color=black,mark=square*] table[x=h,y=EzCP3P2] {Figures/dats/tarabay/FDA/Ez.dat};
          \addplot[color=red,mark=square*] table[x=h,y=EzL] {Figures/dats/tarabay/FDA/Ez.dat};
          \addplot[color=brown,mark=square*] table[x=h,y=EzM] {Figures/dats/tarabay/FDA/Ez.dat};
          \addplot[color=blue,mark=square*] table[x=h,y=EzF] {Figures/dats/tarabay/FDA/Ez.dat};
          \legend{\textsc{M0}-P2P1G1,\textsc{M0}-P2P1G2,\textsc{M0}-P3P2G1,\textsc{M1}-P2P1G1,\textsc{M2}-P2P1G1,\textsc{M3}-P2P1G1}
        \end{semilogyaxis}
      \end{tikzpicture}
    \end{subfigure}
        \begin{subfigure}[b]{.47\linewidth}
      \begin{tikzpicture}[scale=0.85]
        \begin{axis}[
            ylabel near ticks,
            yticklabel pos=right,
            xlabel=z,
            ylabel=$E_q$,
            legend style={at={(0.7,0.)}, anchor=south east}
          ]
          \addplot[color=green,mark=square*] table[x=h,y=EqCP2P1]{Figures/dats/tarabay/FDA/Ez.dat};
          \addplot[color=orange,mark=square*] table[x=h,y=EqCP2P1G2] {Figures/dats/tarabay/FDA/Ez.dat};
          \addplot[color=black,mark=square*] table[x=h,y=EqCP3P2] {Figures/dats/tarabay/FDA/Ez.dat};
          \addplot[color=red,mark=square*] table[x=h,y=EqL] {Figures/dats/tarabay/FDA/Ez.dat};
          \addplot[color=brown,mark=square*] table[x=h,y=EqM] {Figures/dats/tarabay/FDA/Ez.dat};
          \addplot[color=blue,mark=square*]  table[x=h,y=EqF] {Figures/dats/tarabay/FDA/Ez.dat};
          \legend{\textsc{M0}-P2P1G1, \textsc{M0}-P2P1G2,\textsc{M0}-P3P2G1,\textsc{M1}-P2P1G1,\textsc{M2}-P2P1G1,\textsc{M3}-P2P1G1}

        \end{axis}
      \end{tikzpicture}
        \end{subfigure}
        \caption{Validation metrics $E_z$ (left) and $E_Q$ (right) for $\re_t = 500$.}
    \label{fig::results:metrics500}
\end{figure}


\end{document}